\documentclass{article}
 \usepackage{amsmath}

\newtheorem{thm}{Theorem}[section]

\newtheorem{prop}[thm]{Proposition}

\newtheorem{lem}[thm]{Lemma}

\newtheorem{ex}[thm]{Example}

\newcommand{\be}{\begin{equation}}
\newcommand{\ee}{\end{equation}}
\newcommand{\ben}{\begin{enumerate}}
\newcommand{\een}{\end{enumerate}}
\newcommand{\beq}{\begin{eqnarray}}
\newcommand{\eeq}{\end{eqnarray}}
\newcommand{\beqn}{\begin{eqnarray*}}
\newcommand{\eeqn}{\end{eqnarray*}}

\newcommand{\pa}{\partial}

\newcommand{\qed}{\hspace*{\fill}Q.E.D.}  

\begin{document}
\title{On A Class of  Finsler Metrics of
Einstein-Reversibility}
\author{Guojun Yang  }
\date{}
\maketitle
\begin{abstract}
In this paper, we introduce the notion of Einstein-reversibility
for Finsler metrics. We  study a class of $p$-power Finsler
metrics defined by a Riemann metric and 1-form
 which are of Einstein-reversibility. It shows that such a class
 of Finsler metrics of Einstein-reversibility are always Einstein
 metrics. In particular, we show that all $p$-power metrics but Randers metrics, square metrics and
 2-dimensional square-root metrics,
 are always Ricci-flat-parallel. Further,
  the local structure is almost determined for 2-dimensional square-root metrics
 which are Einsteinian (equivalently, of isotropic flag
 curvature), and examples show such metrics are not necessarily
 Ricci-flat.

\


\noindent {\bf MR(2000) subject classification: }
 53B40
\end{abstract}

\section{Introduction}

In Finsler geometry, the Finsler metrics under consideration might
not be reversible. This leads to the
 non-reversibility of geodesics and curvatures. However, for
 certain non-reversible Finsler metrics, the geodesics are
 possibly reversible; if the geodesics are non-reversible, the curvatures
 might be reversible. M. Crampin has shown that a Randers
 metric $F =\alpha+\beta$  has strictly reversible geodesics if and only if
 $\beta$ is parallel (\cite{Cr}). Later, Masca-Sabau-Shimada
 studied $(\alpha,\beta)$-metrics with reversible geodesics (\cite{MSS1}
 \cite{MSS2}). They also study $(\alpha,\beta)$-metrics with
 projectively reversible geodesics. In \cite{ShYa}, the present
 author and Z. Shen introduce a weaker reversibility than geodesic-reversibility and they
   study Randers metrics with reversible Riemann
 curvature and Ricci curvature.

In this paper, we study the reversibility of Einstein scalar in
Finsler geometry. For a Finsler metric $F=F(x,y)$ on a manifold
$M$, the Riemann curvature $R_y : T_xM\to T_xM$ is a family of
linear transformations and the Ricci curvature $Ric(x,y)={\rm
trace} (R_y), \; \forall y\in T_xM$. We can always express the
Ricci curvature $Ric(x,y)$ as $Ric(x,y)=(n-1)\lambda(x,y) F^2$ for
some scalar $\lambda(x,y)$ on $TM$, where $\lambda(x,y)$ is called
the {\it Einstein scalar}. A Finsler metric $F$ is called of {\it
Einstein-reversibility} if the Einstein scalar is reversible,
namely, $\lambda(x,y)=\lambda(x,-y)$. Note that
Ricci-reversibility does not imply  Einstein-reversibility, and
Einstein-reversibility does not imply Ricci-reversibility.
  Clearly,  $F$
is Einstein-reversible  if $F$ is reversible.

If the Einstein scalar $\lambda(x,y)$ is a scalar on $M$, namely,
$\lambda(x,y)=\lambda(x)$, then $F$ is called an {\it Einstein
metric}. Einstein metrics are a natural extension of those in
Riemann geometry and they  have been shown to have similar  good
properties as in Riemann geometry for some special Finsler metrics
(\cite{BR2} \cite{CSZ} \cite{CST} \cite{Shen3} \cite{SYa1}). It is
clear that an Einstein metric is Einstein-reversible. If a Finsler
metric is of scalar flag curvature, then Einstein-reversibility is
equivalent to that the  flag curvature is reversible. In
particular, in two-dimensional case, Einstein-reversibility is
equivalent to the reversibility of the flag curvature. Thus it is
a natural problem to study Finsler metrics which are
Einstein-reversible. This is a difficult problem for general
Finsler metrics. We shall restrict our attention to a class of
$p$-power $(\alpha,\beta)$ Finsler metrics in the following form
 \be\label{ycw01}
 F=\alpha(1+\frac{\beta}{\alpha})^p,
 \ee
 where $p\ne 0$ is a real constant, $\alpha=\sqrt{a_{ij}(x)y^iy^j}$ is
Riemannian and $\beta=b_i(x)y^i$ is a 1-form.

In (\ref{ycw01}), if $p=1$, then $F=\alpha+\beta$  satisfying
$b:=\|\beta\|_{\alpha}<1$ is called a {\it Randers metric}.
Randers metrics with special  curvature properties have been
studied by many people in recent years. If $p=2$, then
$F=(\alpha+\beta)^2/\alpha$ satisfying $b:=\|\beta\|_{\alpha}<1$
is called a {\it square metric}. Square metrics have also been
shown to have some special geometric properties (\cite{CSZ}
\cite{SYa} \cite{SY} \cite{Z}).  If $p=-1$, then
$F=\alpha^2/(\alpha+\beta)$ satisfying $b:=\|\beta\|_{\alpha}<1/2$
is called a {\it Matsumoto metric} introduced by M. Matsumoto in
\cite{Mat}. If $p=1/2$, then $F=\sqrt{\alpha(\alpha+\beta)}$
satisfying $b:=\|\beta\|_{\alpha}<1$ is called a {\it square-root
metric}. We will show below that two-dimensional square-root
metrics have some special properties. The four cases are special
in the proof of Theorem \ref{th1} below.

For a $p$-power $(\alpha,\beta)$-metric in (\ref{ycw01}) which is
of Einstein-reversibility, we have the following theorem.

\begin{thm}\label{th1}
Let $F$ be a $p$-power $(\alpha,\beta)$-metric defined by
(\ref{ycw01}) on an $n$-dimensional manifold.  Then $F$ is
Einstein-reversible if and only if $F$ is an Einstein metric.
Further, except for $p=1,2$, and $p=1/2$ but $n=2$,  $F$ is
Ricci-flat-parallel ($\alpha$ is Ricci-flat and $\beta$ is
parallel with respect to $\alpha$).
\end{thm}

It is shown in \cite{ShYa} that the Ricci curvature of  a Randers
metric is  reversible if and only if the Ricci curvature is
quadratic. Theorem \ref{th1} shows a similar property for
$p$-power metrics of Einstein-reversibility.

An Einstein Randers metric has the similar properties as an
Einstein Riemann metric. Bao-Robles prove that any $n(\ge
3)$-dimensional Einstein Randers metric is of Ricci constant and
in particular is of constant flag curvature when $n=3$
(\cite{BR2}). These results follow from the navigation problem
which shows that an Einstein Randers metric is constructed from an
Einstein Riemann metric and a conformal 1-form with constant
factor.

If a square metric  is an Einstein metric, then it can be verified
directly that $\beta$ is closed. This fact is also shown in Lemma
\ref{lem42} below. Then the local structure of Einstein square
metrics can be determined up to the local structure of Einstein
Riemann metrics (see \cite{CSZ}, therein the Douglas condition can
be cancelled).

In \cite{YZ}, it gives an incomplete characterization for Einstein
Matsumoto metrics. By Theorem \ref{th1}, a Matsumoto metric of
Einstein-reversibility is trivial, namely, Ricci-flat-parallel.
The  discussion in the following  sections says that, among
$p$-power metrics of Einstein-reversibility, $\beta$ is always
closed excluding  Randers metrics and  two-dimensional square-root
metrics.

The local structure can be almost determined for two-dimensional
square-root metrics which are Einstein-reversible. We have the
following theorem.

\begin{thm}\label{th2}
 Let $F=\sqrt{\alpha(\alpha+\beta)}$ be a two-dimensional
 square-root metric which is Einsteinian (equivalently, of
 isotropic flag curvature). Then $\alpha$ and $\beta$ can be
 locally
 determined by
 \be\label{ygj001}
  \alpha=\frac{\sqrt{B}}{(1-B)^{\frac{3}{4}}}\sqrt{\frac{(y^1)^2+(y^2)^2}{u^2+v^2}},\
  \ \ \ \
  \beta=\frac{B}{(1-B)^{\frac{3}{4}}}\frac{uy^1+vy^2}{u^2+v^2},
  \ee
  where
    $0<B=B(x)<1,u=u(x),v=v(x)$ are some scalar functions
   which satisfy the following PDEs:
   \be\label{ygj002}
  u_1=v_2,\ \ u_2=-v_1,\ \ uB_1+vB_2=0,
   \ee
where $u_i:=u_{x^i},v_i:=v_{x^i}$ and $B_i:=B_{x^i}$. Further, the
isotropic flag curvature ${\bf K}$ is given by
 \be\label{ygj003}
{\bf K}=-\frac{(u^2+v^2)\sqrt{1-B}}{2B^2}(B_{11}+B_{22})
-\frac{(u^2+v^2)^2(3B-2)}{4B^3\sqrt{1-B}}\big(\frac{B_1}{v}\big)^2,
 \ee
 where $B_{ij}:=B_{x^ix^j}$.
\end{thm}

We will use a basic result in \cite{Y} \cite{Y1} to prove Theorem
\ref{th2}. It can be shown in Theorem \ref{th2} that $B=contant$
if and only if $\beta$ is closed, if and only if $\alpha$ is flat
and $\beta$ is parallel with respect to $\alpha$ (by
(\ref{ycw084}) below). One can try to find some conditions such
that ${\bf K}=constant$ in (\ref{ygj003}). We have not solved it.

Choose $u=-x^2,\ v=x^1$ and $B=(x^1)^2+(x^2)^2$ in Theorem
\ref{th2}. It is easily shown that $u,v,B$ satisfy (\ref{ygj002}).
We have the following example, which can also be directly
verified.

\begin{ex}\label{ex1}
Let $F=\sqrt{\alpha(\alpha+\beta)}$ be a two-dimensional
 square-root metric, where $\alpha$ and $\beta$ are determined by
 (\ref{ygj001}) with $u=-x^2,\ v=x^1$ and $B=(x^1)^2+(x^2)^2$. Then
 $F$ is of isotropic flag curvature ${\bf K}={\bf K}(x)$ and it
 follows from (\ref{ygj003}) that
  $$
{\bf K}=-\frac{1}{\sqrt{1-(x^1)^2-(x^2)^2}}.
  $$

\end{ex}

In \cite{CST}, it proves that if $\phi(s)$ is a polynomial of $s$
with degree $\ge 2$, then an $n(\ge 2)$-dimensional metric
$F=\alpha\phi(\beta/\alpha)$ is Ricci-flat if $F$ is an Einstein
metric. By taking a limit on the degree, this fact seems still
true for any analytic function $\phi(s)$ (not polynomial).
However, it is not true. By Theorem \ref{th2} and Example
\ref{ex1}, a two-dimensional Einstein square-root metric is
generally not Ricci-flat.

\section{Preliminaries}

For a Finsler metric $F$, the Riemann curvature $R_y=R^i_{\
k}(y)\frac{\pa}{\pa x^i}\otimes dx^k$ is defined by
 \be\label{cw2}
 R^i_{\ k}:=2\frac{\pa G^i}{\pa x^k}-y^j\frac{\pa^2G^i}{\pa x^j\pa
 y^k}+2G^j\frac{\pa^2G^i}{\pa y^j\pa y^k}-\frac{\pa G^i}{\pa y^j}\frac{\pa G^j}{\pa
 y^k},
 \ee
 where $G^i$ are called the geodesic coefficients as follows
 \be\label{cw3}
  G^i:=\frac{1}{4}g^{il}\big\{[F^2]_{x^ky^l}y^k-[F^2]_{x^l}\big\}.
  \ee
Then the Ricci curvature ${\bf Ric}$ is defined by ${\bf
Ric}:=R^k_{\ k}$.

 For a Finsler metric $F$, let
 \be\label{cw6}
 f(x,y):=\frac{R^m_{\ m}}{F^2}=\frac{{\bf Ric}}{F^2}.
 \ee
Then by the definition of Einstein-reversibility, we have
$f(x,y)=f(x,-y)$.

\

To show  the Ricci curvatures of  $p$-power metrics, we need some
notations. For a pair of $\alpha$ and $\beta$, we define the
following quantities:
 $$r_{ij}:=\frac{1}{2}(b_{i|j}+b_{j|i}),\ \ s_{ij}:=\frac{1}{2}(b_{i|j}-b_{j|i}),
 $$
 $$r^i_{\ j}:=a^{ik}r_{kj},\ \  s^i_{\ j}:=a^{ik}s_{kj},\ \ q_{ij}:=r_{im}s^m_{\ j}, \ \ t_{ij}:=s_{im}s^m_{\ j},\ \
 $$
 $$
 r_j:=b^ir_{ij},\ \  s_j:=b^is_{ij},\ \  q_j:=b^iq_{ij}, \ \ r_j:=b^ir_{ij},\ \ t_j:=b^it_{ij},
 $$
 where $b^i:=a^{ij}b_j$, $(a^{ij})$ is the inverse of
 $(a_{ij})$, and $\nabla \beta = b_{i|j} y^i dx^j$  denotes the covariant
derivatives of $\beta$ with respect to $\alpha$.

Here are some of our conventions in the whole paper. For a general
tensor $T_{ij}$ as an example, we define $T_{i0}:=T_{ij}y^j$ and
$T_{00}:=T_{ij}y^iy^j$, etc. We use $a_{ij}$ to raise or lower the
indices of a tensor.

\begin{lem} (\cite{BR2})\label{lem0}
By Ricci identities we have
 $$s_{ij|k}=r_{ik|j}-r_{jk|i}-b^l\bar{R}_{klij},$$
 $$s^k_{\ 0|k}=r^k_{\ k|0}-r^k_{\ 0|k}+b^l\bar{R}ic_{l0},$$
 $$b^ks_{0|k}=r_ks^k_{\ 0}-t_0+b^kb^lr_{kl|0}-b^kb^lr_{k0|l},$$
 $$s^k_{\ |k}=r^k_{\ |k}-t^k_{\ k}-r^i_{\ j}r^j_{\ i}-b^ir^k_{\
 k|i}-b^kb^i\bar{R}ic_{ik}.$$
 where $\bar{R}$ denotes the Riemann curvature tensor of $\alpha$.
\end{lem}

Consider an $(\alpha,\beta)$-metric $F =\alpha \phi (s)$,
$s=\beta/\alpha$. By (\ref{cw3}), the spray coefficients $G^i$ of
$F$ are given by:
  \be\label{ycw02}
  G^i=G^i_{\alpha}+\alpha Q s^i_0+\alpha^{-1}\Theta (-2\alpha Q
  s_0+r_{00})y^i+\Psi (-2\alpha Q s_0+r_{00})b^i,
  \ee
where
 $$
  Q:=\frac{\phi'}{\phi-s\phi'},\ \
  \Theta:=\frac{Q-sQ'}{2\Delta},\ \
  \Psi:=\frac{Q'}{2\Delta},\ \ \Delta:=1+sQ+(b^2-s^2)Q'.
 $$
 By (\ref{cw2}) and (\ref{ycw02}), it gives in \cite{CST}
 \cite{Z} the very long  Ricci curvature {\bf Ric}
 of an $(\alpha,\beta)$-metric. Now for a
 $p$-power metric in (\ref{ycw01}), we have $\phi(s)=(1+s)^p$, and
 then we obtain the corresponding Ricci curvature {\bf Ric}. But
 in general case, the expression is still long. For the simplest
 case of Randers metric, we have

\begin{lem} {\rm(\cite{BR})}\label{lem2.2}
 Let $F=\alpha+\beta$ be an $n$-dimensional Randers
metric on a manifold $M$. Then the Ricci curvature of $F$ is given
by
 \beq\label{cw4}
 {\bf Ric}&=&\widetilde{{\bf Ric}}+2\alpha s^m_{\
 0|m}-2t_{00}-\alpha^2 t^m_{\ m}+\nonumber\\
 &&(n-1)\Big\{\frac{3(r_{00}-2\alpha
 s_0)^2}{4F^2}+\frac{4\alpha(q_{00}-\alpha t_0)-r_{00|0}+2\alpha
 s_{0|0}}{2F}\Big\},
  \eeq
  where $\widetilde{{\bf Ric}}:=\widetilde{R}^m_{\ m}$ denotes the
  Ricci curvature of $\alpha$.
\end{lem}

\section{Randers metrics}

In this section, we will prove Theorem \ref{th1} when
$F=\alpha+\beta$ is a Randers metric.

Now for $F=\alpha+\beta$  in (\ref{cw6}), it follows from
(\ref{cw4}) and $f(x,y)=f(x,-y)$ that
 \beq\label{cw7}
 0&=&4[t^k_{\ k}\beta+s^k_{\
 0|k}-(n-1)t_0](\alpha^2-\beta^2)^3+\Big\{4t^k_{\ k}\beta^3+4\big[2s^k_{\
 0|k}-5(n-1)t_0\big]\beta^2\nonumber\\
 &&-2\big[2Ric_{\alpha}-4t_{00}+3(n-1)(2q_{00}+4s_0^2+s_{0|0})\big]\beta-(n-1)(r_{00|0}+6r_{00}s_0)\Big\}\times\nonumber\\
 &&(\alpha^2-\beta^2)^2-(n-1)\beta\big[16t_0\beta^3+8(s_{0|0}+2q_{00}+9s_0^2)\beta^2+4(r_{00|0}+12r_{00}s_0)\beta\nonumber\\
 &&+6r_{00}^2\big](\alpha^2-\beta^2)-12(n-1)\beta^3(r_{00}+2\beta
 s_0)^2.
 \eeq
By (\ref{cw7}), there is a scalar function $c=c(x)$ such that
 \be\label{cw8}
 r_{00}=-2\beta s_0+2c(\alpha^2-\beta^2).
 \ee
 By (\ref{cw8}) we get
 \beq
 r_{00|0}&=&4\beta s_0^2-2\beta
 s_{0|0}+8c\beta^2s_0+2(c_0-2cs_0-4\beta c^2)(\alpha^2-\beta^2),\label{cw91}\\
 q_{00}&=&-s_0^2-\beta t_0-2c\beta s_0,\label{cw10}
 \eeq
where $c_i:=c_{x^i}$. Now plug (\ref{cw8})--(\ref{cw10}) into
(\ref{cw7}) and then we obtain
 \beq\label{cw11}
 0&=&\big[2s^k_{\
 0|k}-(n-1)(2t_0+4cs_0+3c_0)\big]\beta^2-2\big[Ric_{\alpha}-2t_{00}-\alpha^2t^k_{\
 k}+\nonumber\\
 &&(n-1)(s_{0|0}+s_0^2+4c^2\alpha^2)\big]\beta+\big[2s^k_{\
 0|k}-(n-1)(2t_0+4cs_0+c_0)\big]\alpha^2.
 \eeq
It follows from (\ref{cw11}) that
 \be\label{cw12}
 s^k_{\ 0|k}=\frac{1}{2}(n-1)(2\sigma\beta+2t_0+4cs_0+c_0),
 \ee
where $\sigma=\sigma(x)$ is a scalar function. Substituting
(\ref{cw12}) into (\ref{cw11}) yields
 \be\label{cw13}
 Ric_{\alpha}=2t_{00}+t^k_{\
 k}\alpha^2+(n-1)\big[\sigma(\alpha^2+\beta^2)-4c^2\alpha^2-s_{0|0}-s_0^2-c_0\beta\big].
 \ee

 Finally, plug (\ref{cw8})--(\ref{cw10}), (\ref{cw12})  and (\ref{cw13})
 into the Ricci curvature (\ref{cw4}), and then we obtain
  \be\label{cw14}
 Ric=(n-1)(\sigma-c^2)F^2.
  \ee
Now (\ref{cw14}) shows that $F=\alpha+\beta$ is an Einstein
metric. Therefore, a Randers metric of Einstein-reversibility is
 an Einstein metric.         \qed

\section{Square metrics}\label{sec4}

In this section, we will prove Theorem \ref{th1} when
$F=(\alpha+\beta)^2/\alpha$ is a square metric.

\begin{lem}\label{lem41}
An $n$-dimensional square metric $F=(\alpha+\beta)^2/\alpha$ is
Einstein-reversible if and only if the following hold
 \beq
 Ric_{\alpha}\hspace{-0.6cm}&&=-c^2\Big\{\big[2(2n-5)b^2+5(n-1)\big]\alpha^2-6(n-2)\beta^2\Big\}+\frac{t^k_{\ k}}{2(n-1)^2}\Big\{\big[(n+1)(5n+3)\nonumber\\
 &&+8(n+1)b^2-\frac{(n-1)(n-3)}{b^2}\big]\alpha^2+(n+1)(9n-17)\big(1+\frac{1}{b^2}\big)\beta^2\Big\},\label{cw15}\\
 r_{00}\hspace{-0.6cm}&&=c\big[(1+2b^2)\alpha^2-3\beta^2\big], \  \ t_{00}=\frac{t^k_{\ k}}{(n-1)b^2}(b^2\alpha^2-\beta^2),
  \ s^k_{\ 0|k}=-\frac{t^k_{\ k}}{b^2}\beta,    \ s_0=0,\label{cw16}\\
 c_i\hspace{-0.6cm}&&=-2\Big\{\frac{(n+1)(1+b^2)}{(n-1)^2b^2}t^k_{\ k}+c^2\Big\}b_i,\ \ \ (c_i:=c_{x^i}),\label{cw17}
 \eeq
where $c=c(x)$ is a scalar function.
\end{lem}

{\it Proof} : For the square metric $F=(\alpha+\beta)^2/\alpha$,
let
 $f(x,y)$ defined by (\ref{cw6}).
Then by the definition of Einstein-reversibility, we have
$f(x,y)=f(x,-y)$, which can be firstly written in the following
form
 \be\label{cw18}
 (\cdots)\big[(1+2b^2)\alpha^2-3\beta^2\big]-108\beta^2(\alpha^2-\beta^2)^2(2\alpha^2s_0+\beta
 r_{00})\big[4\alpha^2\beta s_0+(\alpha^2+\beta^2)r_{00}\big]=0,
 \ee
where the omitted term in the bracket is a polynomial in $(y^i)$.
By (\ref{cw18}) we know that $(1+2b^2)\alpha^2-3\beta^2$ is
divided by
 $$
2\alpha^2s_0+\beta r_{00}, \ \ \ \text{or} \ \ \ 4\alpha^2\beta
s_0+(\alpha^2+\beta^2)r_{00}.
 $$

(i). If $(1+2b^2)\alpha^2-3\beta^2$ is divided by
$2\alpha^2s_0+\beta r_{00}$, then
 \be\label{cw19}
2\alpha^2s_0+\beta r_{00}=\theta
\big[(1+2b^2)\alpha^2-3\beta^2\big],
 \ee
where $\theta$ is a 1-form. Clearly, (\ref{cw19}) is equivalent to
 $$
 (2s_0-\theta-2b^2\theta)\alpha^2+\beta(r_{00}+3\beta\theta)=0,
 $$
which implies
 $$r_{00}=-3\beta\theta+(1+2b^2)c\alpha^2,\ \ \
 \theta=c\beta+\frac{2s_0}{1+2b^2},$$
 where $c=c(x)$ is a scalar function.
 In the above, plugging $\theta$ into $r_{00}$ we obtain
  \be\label{cw20}
 r_{00}=c\big[(1+2b^2)\alpha^2-3\beta^2\big]-\frac{6\beta
 s_0}{1+2b^2}.
  \ee

(ii). If $(1+2b^2)\alpha^2-3\beta^2$ is divided by $4\alpha^2\beta
s_0+(\alpha^2+\beta^2)r_{00}$, then
 $$
 4\alpha^2\beta
s_0+(\alpha^2+\beta^2)r_{00}=\theta
\big[(1+2b^2)\alpha^2-3\beta^2\big],
 $$
where $\theta$ is a homogeneous polynomial in $(y^i)$ of degree
two.  Similarly as the analysis in (i) we obtain
 \be\label{cw21}
 r_{00}=c\big[(1+2b^2)\alpha^2-3\beta^2\big]-\frac{6\beta
 s_0}{2+b^2}.
  \ee

Starting from (\ref{cw20}) or (\ref{cw21}), the discussions are
similar. We only consider one case. In the following, we assume
(\ref{cw21}) holds.

By (\ref{cw21}) we can obtain the expressions of the following
quantities:
 $$
r_{00},\ r^k_{\ k},\ r, \ r_{00|0}, \ r_{0|0}, \ r_{00|k},
 \ r_{k|0}, \ q_{ij}, \ q_k,\ b^mq_{0m},
\ etc.
 $$
 For example we have
 $$
 q_{ij}=c\big[(1+2b^2)s_{ij}-3b_is_j\big]-\frac{3}{2+b^2}(b_it_j+s_is_j).
 $$

Plug all the above quantities into (\ref{cw18}) and then
multiplied by $1/[(1+2b^2)\alpha^2-3\beta^2]$ the equation
(\ref{cw18}) is written as
 \be\label{cw22}
 (\cdots)\big[(1+2b^2)\alpha^2-3\beta^2\big]+(n-4)s_0^2\beta\alpha^6(\alpha^2+\beta^2)(\alpha^2-\beta^2)^3=0,
 \ee
where the omitted term in the bracket is a polynomial in $(y^i)$.
By (\ref{cw22}), it is clear that
 $$s_0=0, \ \ \ \text{or} \ \ \ n=4.$$
We will show that $n=4$ still implies  $s_0=0$.

\

\noindent {\bf Case I:} Assume $s_0=0$. Then by (\ref{cw21}) we
get $r_{00}$ given in (\ref{cw16}). By $s_0=0$ we have
 $$
 s_{0|0}=0,\ s_{0|k}=0,\ t_0=0,\ q_{00}=0,\ q_0=0,\ b^kq_{0k}=0.
 $$
Then by (\ref{cw22}) multiplied by
$1/[(1+2b^2)\alpha^2-3\beta^2]^3$ we obtain
 \be\label{cw23}
 (\cdots)\alpha^2-\beta^3\big[4Ric_{\alpha}+(9n-17)c_0\beta-2(3n-7)c^2\beta^2\big],
 \ee
where the omitted term in the bracket is a polynomial in $(y^i)$.
Clearly, (\ref{cw23}) shows
 \be\label{cw24}
 Ric_{\alpha}=\sigma\alpha^2+\frac{(3n-7)c^2}{2}\beta^2-\frac{9n-17}{4}c_0\beta,
 \ee
where $\sigma=\sigma(x)$ is a scalar function. Plugging
(\ref{cw24}) into (\ref{cw23}) yields
 \be\label{cw25}
 A\alpha^2+B\beta=0,
 \ee
where
 \beqn
 A:\hspace{-0.7cm}&&=\big\{[8(3-2n)b^2-22(n-1)]c^2+4(2t^k_{\
 k}-2b^kc_k-\sigma)\big\}\beta+4s^k_{\ 0|k}-(n-1)c_0,\\
 B:\hspace{-0.7cm}&&=\big\{[8(3-2n)b^2-30(n-1)]c^2-4(2b^kc_k+\sigma)\big\}\beta^2+[12s^k_{\
 0|k}-5(n-1)c_0]\beta+16t_{00}.
 \eeqn
By (\ref{cw25}) we easily get
 \be
 s^k_{\ 0|k}=\Big\{(2b^kc_k-2t^k_{\
 k}-4d+\sigma)+\big[\frac{11}{2}(n-1)+2(2n-3)b^2\big]c^2\Big\}\beta+\frac{n-1}{4}c_0,\label{cw26}
 \ee
 \be
 t_{00}=d\alpha^2+\Big\{\big[\frac{15}{8}(n-1)+\frac{1}{2}(2n-3)b^2\big]c^2+\frac{1}{2}b^kc_k+\frac{1}{4}\sigma\Big\}\beta^2
 +\Big\{\frac{5}{16}(n-1)c_0-\frac{3}{4}s^k_{\
 0|k}\Big\}\beta,\label{cw27}
 \ee
where $d=d(x)$ is a scalar function.

Contracting (\ref{cw27}) by $a^{ij}$ we get
 \be\label{cw28}
 d=\frac{1}{8(n+3b^2)}\Big\{(8b^2+1-n)b^kc_k+2b^2[4(2n-3)b^2+9(n-1)]c^2+4b^2\sigma-4(3b^2-2)t^k_{\
 k}\Big\}.
 \ee
Contracting (\ref{cw27}) by $b^ib^j$ we get
 \be\label{cw29}
 \sigma=\big[2(3-2n)b^2-\frac{9}{2}(n-1)\big]c^2+\frac{n-1-8b^2}{4b^2}b^kc_k+\frac{2+3(n+1)b^2}{(n-1)b^2}t^k_{\
 k}.
 \ee
Plugging (\ref{cw28}) and (\ref{cw29}) into (\ref{cw27}) and then
contracting (\ref{cw27}) by $b^j$ and using $s_0=0$ we obtain
 \be\label{cw30}
 c_i=\frac{b^kc_k}{b^2}b_i.
 \ee

Plug (\ref{cw28}), (\ref{cw29}) and  (\ref{cw30}) into
(\ref{cw27}) and then (\ref{cw27}) is reduced to the $t_{00}$
given in (\ref{cw16}). Plugging (\ref{cw28})--(\ref{cw30}) into
(\ref{cw26}) we have
 \be\label{cw31}
 s^k_{\
 0|k}=\Big\{(n-1)c^2+\frac{n-1}{2b^2}b^kc_k+\frac{2+(n+1)b^2}{(n-1)b^2}t^k_{\
 k}\Big\}\beta.
 \ee
Contracting (\ref{cw31}) by $b^i$ and using $s_0=0$ we get
 \be\label{cw32}
b^kc_k=-2b^2c^2-\frac{2(n+1)(1+b^2)}{(n-1)^2}t^k_{\ k}.
 \ee

Now using (\ref{cw32}), it follows from (\ref{cw31}) that $s^k_{\
0|k}$ is given in (\ref{cw16}), and we get (\ref{cw17}) from
(\ref{cw30}) and (\ref{cw32}). Finally, (\ref{cw15}) follows from
(\ref{cw24}), (\ref{cw29}) and (\ref{cw30}).

\

\noindent {\bf Case II:} Assume $n=4$. In this case, we will prove
$s_0=0$. The idea of proof is similar as that in Case I, but the
computations will be more complicated. For briefness, some details
will be omitted. In the following, all $A_i,B_i$ are homogeneous
polynomials in $(y^i)$.

Plug $n=4$ into (\ref{cw22}) multiplied by
$1/[(1+2b^2)\alpha^2-3\beta^2]$, and then (\ref{cw22}) can be
written as
 \be\label{cw33}
 A_1\alpha^2+B_1\beta^7=0.
 \ee
It follows from (\ref{cw33}) that $\alpha^2|B_1$. So we have
 \be\label{cw34}
 Ric_{\alpha}=\sigma\alpha^2+\frac{5c^2}{2}\beta^2-\Big\{\frac{3(25b^2+12)}{2(2+b^2)^2}cs_0
 +\frac{19}{4}c_0\Big\}\beta-\frac{(56+85b^2)}{(2+b^2)^3}s_0^2-\frac{19}{2(2+b^2)}s_{0|0},
 \ee
where $\sigma=\sigma(x)$ is a scalar function. Plug (\ref{cw34})
into (\ref{cw33}) and then we have
 \be\label{cw35}
A_2\alpha^2+B_2\beta^5=0.
 \ee
Similarly, (\ref{cw35}) shows $\alpha^2|B_2$, and then we have
 \beq\label{cw36}
 t_{00}&=&d\alpha^2+\Big\{\frac{5b^2c^2}{2}+\frac{\sigma}{4}+\frac{45c^2}{8}+\frac{b^kc_k}{2}\Big\}\beta^2
 +\frac{11b^4-16b^2-58}{3(2+b^2)^3}s_0^2+\frac{1+2b^2}{6(2+b^2)}s_{0|0}\nonumber\\
 &&\Big\{\frac{176b^4+425b^2+506}{24(2+b^2)^2}cs_0-\frac{3s^k_{\
 0|k}}{4}+\frac{b^ks_{0|k}}{2+b^2}+\frac{15}{16}c_0+\frac{17}{2(2+b^2)}t_0\Big\}\beta,
 \eeq
where $d=d(x)$ is a scalar function. Contracting (\ref{cw36}) by
$a^{ij}$ we get
 \beq\label{cw37}
 \sigma&=&-\frac{16d}{b^2}-\frac{8b^2+15}{4b^2}b^kc_k-5\big(2b^2+\frac{9}{2}\big)c^2+\frac{t^k_{\
 k}}{b^2}-\frac{(13b^2+20)s^k_{|k}}{3b^2(2+b^2)}\nonumber\\
 &&+\frac{2(5b^4+176b^2+296)s_ks^k}{3b^2(2+b^2)^3}.
 \eeq
Contracting (\ref{cw36}) by $b^ib^j$ and using (\ref{cw37}) we get
 \be\label{cw38}
d=\frac{1}{3}t^k_{\
k}-\frac{1+2b^2}{18(2+b^2)}s^k_{|k}-\frac{4b^6-46b^4-103b^2-26}{9b^2(2+b^2)^3}s_ks^k.
 \ee
Plug (\ref{cw36})--(\ref{cw38}) into (\ref{cw35}) and we get
 \be\label{cw39}
 A_3\alpha^2+B_3\beta^3=0.
 \ee
Similarly, by (\ref{cw39}) we obtain
 \be\label{cw40}
 s_{0|0}=e\alpha^2+\cdots,
 \ee
where $e=e(x)$ is a scalar function. Contracting (\ref{cw40}) by
$b^ib^j$ we can solve $e$. Plugging $e$ into (\ref{cw40}) and then
contracting (\ref{cw40}) by $a^{ij}$ we have
 \be\label{cw41}
 t^k_{\
 k}=\frac{(b^4-26b^2-38)s^k_{|k}}{45(2+b^2)(1+b^2)}-\frac{92b^8+485b^6+648b^4+724b^2+832}{45b^2(2+b^2)^3(1+b^2)}s_ks^k
 -\frac{9b^2c^2}{5(1+b^2)}-\frac{9b^kc_k}{10(1+b^2)}.
 \ee
Plugging (\ref{cw40}), (\ref{cw41}) and $e$ into (\ref{cw39})
yields
 \be\label{cw42}
 A_4\alpha^2+B_4\beta=0.
 \ee
Similarly, by (\ref{cw42}) we can obtain
 \be\label{cw43}
 s_0^2=u\alpha^2+\cdots,
 \ee
where $u=u(x)$ is a scalar function. Contracting (\ref{cw43}) by
$a^{ij}$ we can solve $u$. Then Contracting (\ref{cw43}) by
$b^ib^j$ and  using $u$ and (\ref{cw41}) we obtain
 \be\label{cw44}
 s^k_{|k}=\frac{2(b^6+12b^4+38b^2+24)}{(2+b^2)(b^4+5b^2+6)b^2}s_ks^k.
 \ee
Plug (\ref{cw43}), (\ref{cw44}) and $u$ into (\ref{cw42}), and
then finally we obtain $s^k_{\ 0|k}$. Then contracting $s^k_{\
0|k}$ by $b^i$ and using (\ref{cw41}) and (\ref{cw44}), we obtain
 \be\label{cw45}
 \frac{8(7b^2+20)(b^2-1)^3}{9(3+b^2)(1+2b^2)^2(2+b^2)}s_ks^k=0.
 \ee
Clearly, (\ref{cw45}) shows $s_k=0$.         \qed

\begin{lem}\label{lem42}
Let $F=(\alpha+\beta)^2/\alpha$ be an $n$-dimensional square
metric of Einstein-reversibility. Then $\beta$ is closed.
\end{lem}

{\it Proof :} We use the fourth formula of Lemma \ref{lem0} to
prove our result. By the fourth formula in Lemma \ref{lem0}, we
first compute $r^k_{\ k|i}$, $r^k_{\ m}r^m_{\ k}$, $r^k_{\ |k}$
and $b^kb^l\bar{R}ic_{kl}$. By (\ref{cw15}),  we have
 \be\label{ycw046}
b^kb^l\bar{R}ic_{kl}=\frac{t^m_{\
m}}{2(n-1)}\big[(n+1)b^2(9b^2+14)-n+3\big]+(n-1)c^2b^2(2b^2-5)
 \ee
It follows from $r_{00}$ in (\ref{cw16}) that
 \beq
r^k_{\ |k}&=&(1-b^2)\big[(2nb^2-5b^2+n)c^2+b^mc_m\big],\label{ycw047}\\
 r^k_{\
k|i}&=&(n-3b^2+2nb^2)c_i+2(2n-3)(1-b^2)c^2b_i.\label{ycw048}
 \eeq
Therefore, by (\ref{ycw046})-(\ref{ycw048}), (\ref{cw17}) and
$r_{00}$ in (\ref{cw16}), we obtain from the fourth formula of
Lemma \ref{lem0}
 \be\label{ycw049}
 s^k_{\
|k}=\frac{(n+1)(1-b^2)(3+b^2)}{2(n-1)}t^m_{\ m}.
 \ee
But $s^k_{\ |k}=0$ since $s_0=0$. So (\ref{ycw049}) implies
$t^m_{\ m}=0$. It is easy to show that $\beta$ is closed if and
only if $t^m_{\ m}=0$.          \qed

\begin{prop}
An $n$-dimensional square metric $F=(\alpha+\beta)^2/\alpha$ is
Einstein-reversible if and only if $F$ is Ricci-flat.
\end{prop}

{\it Proof :} By Lemma \ref{lem42}, (\ref{cw15})-(\ref{cw17}) are
reduced to
 \beq
 Ric_{\alpha}\hspace{-0.6cm}&&=-c^2\Big\{\big[2(2n-5)b^2+5(n-1)\big]\alpha^2-6(n-2)\beta^2\Big\},\label{cw015}\\
 b_{i|j}\hspace{-0.6cm}&&=c\big[(1+2b^2)a_{ij}-3b_ib_j\big], \label{cw016}\\
 c_i\hspace{-0.6cm}&&=-2c^2b_i.\label{cw017}
 \eeq
Plug (\ref{cw015})-(\ref{cw017}) into the Ricci curvature of the
square metric $F$ and then we obtain ${\bf Ric}=0$. So $F$ is a
Ricci-flat Einstein metric.

\section{Matsumoto metrics}

In this section, we will prove Theorem \ref{th1} when
$F=\alpha^2/(\alpha+\beta)$ is a Matsumoto metric.

\begin{prop}\label{prop51}
An $n$-dimensional Matsumoto metric $F=\alpha^2/(\alpha+\beta)$ is
Einstein-reversible if and only if $F$ is Ricci-flat-parallel.
\end{prop}

{\it Proof} : For the square metric $F=\alpha^2/(\alpha+\beta)^2$,
let
 $f(x,y)$ defined by (\ref{cw6}).
Then by the definition of Einstein-reversibility, we have
$f(x,y)=f(x,-y)$, which multiplied by
$\alpha^4(4s^2-1)^3[(1+2b^2)^2-9s^2]^4$ can be firstly written in
the following form
 \be\label{ycw053}
 (\cdots)(1+2b^2+3s)+8748s^4(1+s)^4(1+2s)^3(1-2s)^3\big[(1+2s)r_{00}+2\alpha s_0\big]^2=0,
 \ee
where $s:=\beta/\alpha$ and the omitted term in the bracket is in
the following form
 \be\label{ycw00}
  f_0(s)+f_1(s)\alpha+\cdots +f_m(s)\alpha^m,
  \ee
  for some integer $m$ and $f_i(s)$'s being polynomials of $s$ with coefficients being
  homogenous polynomials in $(y^i)$. By (\ref{ycw053}) we have
   $$\alpha\big[(1+2s)r_{00}+2\alpha
   s_0\big]=(g_0+g_1\alpha+g_2\alpha^2)\alpha(1+2b^2+3s),$$
   which is equivalent to
  \be\label{ycw055}
   g_2(1+2b^2)\alpha^3+\big[(1+2b^2)g_1-2s_0+3g_2\beta\big]\alpha^2-\big[r_{00}-(1+2b^2)g_0-3g_1\beta\big]\alpha-(2r_{00}-3g_0)\beta=0,
   \ee
   where $g_0$ is a polynomial of degree two in $(y^i)$, $g_1$ is
   a 1-form and $g_2=g_2(x)$ is a scalar function. Now it easily
   follows from (\ref{ycw055}) that
    \be\label{ycw056}
   r_{00}=c\big[(1+2b^2)^2\alpha^2-9\beta^2\big]+\frac{18}{(1+2b^2)(1-4b^2)}\beta
   s_0,
    \ee
    where $c=c(x)$ is a scalar function.

By (\ref{ycw056}) we can obtain the expressions of the following
quantities:
 $$
r_{00},\ r^k_{\ k},\ r, \ r_{00|0}, \ r_{0|0}, \ r_{00|k},
 \ r_{k|0}, \ q_{ij}, \ q_k,\ b^mq_{0m},
\ etc.
 $$
Plug all the above quantities into (\ref{ycw053}) and then
 the equation (\ref{ycw053}) multiplied by $\alpha^7[(1+2b^2)^2-9s^2]^{-4}$ is written as
 \be\label{ycw057}
 (\cdots)(\alpha^2-4\beta^2)+2048(1-4b^2)^2(1+2b^2)^3\beta^9\big[(4b^2-1)t_{00}+4s_0^2\big]=0,
 \ee
where the omitted term in the bracket is a polynomial in $(y^i)$.
 It follows from (\ref{ycw057}) that
  \be\label{ycw058}
 t_{00}=d(\alpha^2-4\beta^2)+\frac{4}{1-4b^2}s_0^2,
  \ee
where $d=d(x)$ is a scalar function. By (\ref{ycw058}) we get
 \be\label{ycw059}
 t_0=d(1-4b^2)\beta,\ \ \  t^m_{\
 m}=d(n-4b^2)+\frac{4s_ms^m}{1-4b^2},\ \ \  s_ms^m=-d(1-4b^2)b^2.
 \ee
Plugging (\ref{ycw058}) and (\ref{ycw059}) into (\ref{ycw057})
yields
 \be\label{ycw060}
 (\cdots)(\alpha^2-4\beta^2)+1536(n-4)(1-4b^2)(1+2b^2)^3\beta^7\big[d(1+16b^4-8b^2)\beta^2-s_0^2\big]=0,
 \ee
where the omitted term in the bracket is a polynomial in $(y^i)$.
Eq.(\ref{ycw060}) implies
 $$n=4,\ \ \ {\text or} \ \ \ \ \
 d(1+16b^4-8b^2)\beta^2-s_0^2=-\theta(\alpha^2-4\beta^2),$$
 where $\theta=\theta(x)$ is a scalar function. We will show that $\beta$ is
 parallel with respect to $\alpha$ in the above two cases.

\

\noindent{\bf Case I:} Assume $n=4$. In this case, (\ref{ycw060})
can be written as
 \be\label{ycw061}
 (\cdots)(1-2s)-\alpha(1-4b^2)^2(1+2b^2)^3\big[(1-4b^2)(s^k_{\
 0|k}+8d\alpha)+8c(1+6b^2+2b^4)s_0+4b^ks_{0|k}\big],
 \ee
where the omitted term in the bracket is in the form
(\ref{ycw00}). By (\ref{ycw061}) we have
 \be\label{ycw062}
(1-4b^2)(s^k_{\
 0|k}+8d\alpha)+8c(1+6b^2+2b^4)s_0+4b^ks_{0|k}=\theta\alpha(1-2s),
 \ee
where $\theta=\theta(x)$ is a scalar function. It follows from
(\ref{ycw062}) that
 \be\label{ycw063}
 s^k_{\
 0|k}=-16d\beta-\frac{4}{1-4b^2}\big[2c(1+6b^2+2b^4)s_0+b^ks_{0|k}\big].
 \ee
 Plug (\ref{ycw063}) into (\ref{ycw061}) and then (\ref{ycw061})
 has the following form
  \be\label{ycw064}
 A\alpha^4+B\alpha^2+C=0,
  \ee
where $A,B,C$ are polynomials in $(y^i)$. By $\alpha^2|C$ we get
 \beq\label{ycw065}
 s_{0|0}&=&e\alpha^2+\frac{6(1-4b^2)(1+2b^2)(5b^2-17)c^2\beta^2}{5}+\Big\{\frac{(1+2b^2)(1-4b^2)c_0}{2}\nonumber\\
 &&-\frac{(80b^4-466b^2-343)cs_0}{5(1+2b^2)}\Big\}\beta+\frac{2(160b^6+672b^4+228b^2-169)s_0^2}{5(1+2b^2)^2(1-4b^2)^2},
 \eeq
 where $e=e(x)$ is a scalar function.
Plug (\ref{ycw065}) into (\ref{ycw064}) and then similarly we
obtain
 \beq\label{ycw066}
 Ric_{\alpha}&=&f\alpha^2-\frac{9(352b^2-427)s_0^2}{5(1+2b^2)^2(1-4b^2)^2}+\frac{6(200b^6+226b^4-1223b^2-688)c\beta
 s_0}{5(1+2b^2)^2(1-4b^2)}\nonumber\\
 &&+\Big\{3b^kc_k+(\frac{423b^2}{5}+\frac{2259}{10})c^2-\frac{21d}{1+2b^2}-\frac{15e}{(1+2b^2)(1-4b^2)}\Big\}\beta^2\nonumber\\
 &&-\Big\{\frac{6b^ks_{0|k}}{(1+2b^2)(1-4b^2)}+3(b^2+\frac{1}{4})c_0\Big\}\beta,
 \eeq
where $f=f(x)$ is a scalar function, and
 \beq\label{ycw067}
 b^ks_{0|k}&=&\frac{(80b^6-92b^4-56b^2-13)cs_0}{2(1+2b^2)}-\frac{(1+2b^2)(1+4b^2)(1-4b^2)c_0}{8}\nonumber\\
 &&-\Big\{(1+2b^2)(1-4b^2)\big[\frac{(64b^6+20b^4-10b^2-11)c^2}{4}+\frac{f}{6}-\frac{1-4b^2}{6}b^kc_k\big]\nonumber\\
 &&+\frac{(1-4b^2)(8b^2+19)d}{6}+\frac{8}{3}(1-B)e\Big\}.
 \eeq
Using $b^mb^ks_{m|k}=-r^ks_k+s_ks^k$, (\ref{ycw056}) and
(\ref{ycw059}), it follows from (\ref{ycw067}) that
 \beq\label{ycw068}
 f&=&(\frac{33}{2}+15b^2-30b^4-96b^6)c^2+\frac{(64b^6+120b^4-36b^2-13)d}{(1+2b^2)^2(1-4b^2)}\nonumber\\
 &&-\frac{16(1-b^2)e}{(1+2b^2)(1-4b^2)}-\frac{16b^4+8b^2+3}{4b^2}b^kc_k.
 \eeq
Similarly, it follows from (\ref{ycw065}) that
 \be\label{ycw069}
 e=-(1-4b^2)(1+2b^2)\Big\{\frac{6(5b^2-17)b^2c^2}{5}+\frac{1}{2}b^kc_k\Big\}+\frac{d(8b^4+11b^2-1)}{1+2b^2}.
 \ee
Further by (\ref{ycw065}) we have
 \beq\label{ycw070}
 s^m_{|m}&=&4e+(1+2b^2)(1-4b^2)\Big\{\frac{b^kc_k}{2}+\frac{6(5b^2-17)b^2c^2}{5}\Big\}\nonumber\\
 &&-\frac{2d(1-4b^2)b^2(160b^6+672b^4+228b^2-169)}{5(1+2b^2)^2(1-4b^2)^2}.
 \eeq
Using $b^ks^m_{\ k|m}=-s^m_{|m}-t^m_{\ m}$,
(\ref{ycw068})-(\ref{ycw070}) and (\ref{ycw059}), it follows from
(\ref{ycw063}) that
 \be\label{ycw071}
 d=-\frac{(1+2b^2)^3(1-4b^2)^2}{484b^4+166b^2-83}\big\{3(5b^2-17)c^2+\frac{5}{4b^2}b^kc_k\big\}.
 \ee

Now plug (\ref{ycw056}), (\ref{ycw059}), (\ref{ycw066}),
(\ref{ycw068})-(\ref{ycw071}) into the fourth formula in Lemma
\ref{lem0}, and then we obtain
 \beq\label{ycw072}
 &&6b^2(1440b^{10}+29316b^8+55102b^6+1446b^4-1791b^2-463)c^2\nonumber\\
 &=&-(720b^8+4522b^6+4875b^4-681b^2-364)b^mc_m.
 \eeq
Plugging (\ref{ycw056}), (\ref{ycw059}),
(\ref{ycw068})-(\ref{ycw071}) into the third formula in Lemma
\ref{lem0} yields
 \be\label{ycw073}
 12b^2(48b^4+10b^2+5)cs_0+(1-4b^2)(1+2b^2)(1+20b^2)(b^2c_0-b^mc_m\beta)=0.
 \ee
Similarly, plugging (\ref{ycw056}), (\ref{ycw066}) and
(\ref{ycw063})  into the second formula in Lemma \ref{lem0}, and
using (\ref{ycw059}) and (\ref{ycw067})-(\ref{ycw071}), we get an
equation, and from this equation and (\ref{ycw072}) we eliminate
the terms including $c^2$ and then the final result gives
 \beq\label{ycw074}
 &&3(5120b^8-6848b^6+6832b^4+1091b^2+150)b^2cs_0\nonumber\\
 &=&10(1-4b^2)(32b^4+16b^2+3)(1+2b^2)^2(b^2c_0-b^mc_m\beta).
 \eeq
Now it easily follows from (\ref{ycw073}) and (\ref{ycw074}) that
 \be\label{ycw075}
 3b^2(1+8b^2)(28160b^8-1440b^6+26484b^4+3691b^2+750)cs_0=0.
 \ee
Since $b^2<1/4$, (\ref{ycw075}) implies $c=0$ or $s_0=0$. If
$c=0$, then by (\ref{ycw071}) we have $d=0$, and thus by
(\ref{ycw059}) we have $s_m=0$ and $t^m_{\ m}=0$. Then $r_{00}=0$
by (\ref{ycw056}) and $t^m_{\ m}=0$ imply that $\beta$ is parallel
with respect to $\alpha$. If $s_0=0$, then $d=0$ by
(\ref{ycw059}), and now it easily follows from (\ref{ycw071}) and
 (\ref{ycw072}) that $c=0$. Therefore, $\beta$ is also parallel
with respect to $\alpha$. Further, by (\ref{ycw066}) we have
$Ric_{\alpha}=0$.

\

\noindent{\bf Case II:} Assume $n\ne 4$. In this case, by
(\ref{ycw060}) we have
 \be\label{ycw076}
 d(1+16b^4-8b^2)\beta^2-s_0^2=-\theta(\alpha^2-4\beta^2),
 \ee
 for some scalar function $\theta=\theta(x)$. By contractions on
 (\ref{ycw076}), we have
  \be\label{ycw077}
 \theta=-db^2(1-4b^2), \ \ \ (n-2)db^2(1-4b^2)=0.
  \ee

\noindent{\bf Case IIA:} Assume $n=2$. Plugging (\ref{ycw077})
into (\ref{ycw076}) gives
 $$
 s_0^2=-d(1-4b^2)(b^2\alpha^2-\beta^2).
 $$
Since $n=2$, we can always put
 $$
 Ric_{\alpha}=e\alpha^2,
 $$
for some scalar function $e=e(x)$. Then similar to the discussion
in Case I, we can show that $F$ is Ricci-flat-parallel. The
details are omitted. But the computations show that this case is
much simpler.

\

\noindent{\bf Case IIB:} Assume $n\ne 2$. By (\ref{ycw077}) we
have $d=0$. Then (\ref{ycw059}) shows $\beta$ is closed. So in
this case, (\ref{ycw060}) is reduced to the following form
 \be\label{ycw078}
 A\alpha^4+B\alpha^2+C=0,
  \ee
where $A,B,C$ are polynomials in $(y^i)$. By $\alpha^2|C$ we get
 \be\label{ycw079}
 c_0=-\frac{6c^2(18-6b^2-13n+4nb^2)}{2n-3}\beta.
 \ee
Plug (\ref{ycw079}) into  (\ref{ycw078}) and then (\ref{ycw078})
becomes
 \beq\label{ycw080}
 Ric_{\alpha}&=&-\Big\{\frac{6(22n^2-13n-18)b^4+6(5n^2+10n-24)b^2+(n-1)(17n-21)}{2(2n-3)}+\nonumber\\
 &&16(n-1)b^6\Big\}c^2\alpha^2+\frac{3\big[2(7n^2+20n-36)b^2+67n^2-112n+27\big]}{2(2n-3)}c^2\beta^2.
 \eeq
Now plugging (\ref{ycw056}) into  the third formula in Lemma
\ref{lem0} and using the fact that $\beta$ is closed, we obtain
 \be\label{ycw081}
 c_0=\frac{b^mc_m}{b^2}\beta.
 \ee
By (\ref{ycw079}) and (\ref{ycw081}) we have
 \be\label{ycw082}
 6b^2(4nb^2-6b^2-13n+18)c^2+(2n-3)b^mc_m=0.
 \ee
 Plugging (\ref{ycw056}) and (\ref{ycw080}) into  the fourth formula in Lemma
\ref{lem0} and using the fact that $\beta$ is closed, we obtain
 \be\label{ycw083}
 3b^2\big[(64n-96)b^6-(46n-60)b^4+(73n-99)b^2+17n-27\big]c^2+2(2n-3)(1+2b^2)^2b^mc_m=0.
 \ee
It easily follows from (\ref{ycw082}) and (\ref{ycw083}) that
 $$
 3b^2\big[(98n-132)b^4+(265n-363)b^2+69n-99\big]c^2=0,
 $$
 which implies $c=0$. Thus by (\ref{ycw056}), it is obvious that
 $\beta$ is parallel with respect to $\alpha$. Eq.(\ref{ycw080})
 shows that $Ric_{\alpha}=0$.          \qed

\section{Square-root metrics}

In this section, we will prove Theorem \ref{th1} when
$F=\sqrt{\alpha(\alpha+\beta)}$ is a square-root metric.

\begin{prop}\label{prop61}
An $n$-dimensional square-root metric
$F=\sqrt{\alpha(\alpha+\beta)}$ is Einstein-reversible if and only
if one of the following cases holds

 \ben
  \item[{\rm (i)}] ($n\ge 3$) $F$ is Ricci-flat-parallel.

  \item[{\rm (i)}] ($n=2$) $F$ satisfies
    \beq
     r_{00}&=&\frac{6}{b^2-4}\beta s_0,\label{ycw084}\\
     t_{00}&=&4\theta(4\alpha^2-\beta^2)-\frac{s_0^2}{b^2-4},\ \ \ \ \  s_0^2=4(b^2-4)\theta(b^2\alpha^2-\beta^2),\label{ycw085}\\
    s^k_{\ 0|k}&=&-\frac{32(b^4+7b^2-26)\theta+(b^2-4)^2\lambda}{2(2+b^2)(b^2-4)}\beta+\frac{b^ks_{0|k}}{b^2-4},\label{ycw086}\\
     s_{0|0}&=&\frac{2\big[(5b^6+12b^4-60b^2+16)\alpha^2-6(b^4+2b^2-12)\beta^2\big]s_ks^k}{b^2(2+b^2)(b^2-4)^2}+\nonumber\\
    && \frac{\lambda(b^2-4)(b^2\alpha^2-\beta^2)}{2(2+b^2)},\label{ycw087}
    \eeq
    where $\theta=\theta(x)$ is a scalar function and $\lambda=\lambda(x)$ is the sectional curvature of $\alpha$.
    In this case, $F$ has isotropic flag curvature ${\bf K}={\bf K}(x)$ given by
     \be\label{ycw088}
       {\bf K}=\frac{2(\lambda-32\theta)}{2+b^2}.
     \ee
    We will show (\ref{ycw085})-(\ref{ycw087}) automatically hold
    in the last section.
 \een
\end{prop}

{\it Proof} : For the square-root metric
$F=\sqrt{\alpha(\alpha+\beta)}$, let
 $f(x,y)$ defined by (\ref{cw6}).
Then by the definition of Einstein-reversibility, we have
$f(x,y)=f(x,-y)$, which multiplied by
$\alpha^4(1-s)(s^2-4)^3[(4-b^2+3s^2)^2-36s^2]^4$ can be firstly
written in the following form
 \be\label{ycw089}
 (\cdots)(4-b^2+6s+3s^2)+s^4(1-s)(1+s)^4(2+s)^3(2-s)^3\big[(2+s)r_{00}-2\alpha s_0\big]^2=0,
 \ee
where $s:=\beta/\alpha$ and the omitted term in the bracket is in
a similar form as in (\ref{ycw00}). It follows from (\ref{ycw089})
that
 \be\label{ycw090}
 \alpha\big[(2+s)r_{00}-2\alpha
 s_0\big]=(g_0+g_1\alpha)\alpha^2(4-b^2+6s+3s^2),
 \ee
where $g_0$ is a 1-form and $g_1=g_1(x)$ is a scalar function. By
(\ref{ycw090}) we can easily obtain $r_{00}$ given by
(\ref{ycw084}).

 By (\ref{ycw084}) we can obtain the expressions of
the following quantities:
 $$
r_{00},\ r^k_{\ k},\ r, \ r_{00|0}, \ r_{0|0}, \ r_{00|k},
 \ r_{k|0}, \ q_{ij}, \ q_k,\ b^mq_{0m},
\ etc.
 $$
Plug all the above quantities into (\ref{ycw089}) and then
 the equation (\ref{ycw089}) multiplied by $(b^2-4)^3\alpha^{-2}[(4-b^2+3s^2)^2-36s^2]^{-4}$ is written as
 \be\label{ycw091}
 (\cdots)(4\alpha^2-\beta^2)+(b^2-4)^2\beta^7\big[(b^2-4)t_{00}+s_0^2\big]=0,
 \ee
where the omitted term in the bracket is a polynomial in $(y^i)$.
 It follows from (\ref{ycw091}) that $t_{00}$ in (\ref{ycw085})
 holds for some scalar function $\theta=\theta(x)$. By $t_{00}$ in (\ref{ycw085}) we
get
 \be\label{ycw093}
 t_0=-4\theta(b^2-4)\beta,\ \ \  t^m_{\
 m}=-4\theta(b^2-4n)-\frac{s_ms^m}{b^2-4},\ \ \  s_ms^m=4\theta(b^2-4)b^2.
 \ee
Plugging $t_{00}$ in (\ref{ycw085}) and (\ref{ycw093}) into
(\ref{ycw091}) yields
 \be\label{ycw094}
 (\cdots)(4\alpha^2-\beta^2)+(n-4)(b^2-4)\beta^5\big[\theta(b^2-4)^2\beta^2-s_0^2\big]=0,
 \ee
where the omitted term in the bracket is a polynomial in $(y^i)$.
Eq.(\ref{ycw094}) implies
 $$n=4,\ \ \ {\text or} \ \ \ \ \
 \theta(b^2-4)^2\beta^2-s_0^2=c(4\alpha^2-\beta^2),$$
 where $c=c(x)$ is a scalar function.

\

\noindent{\bf Case I:} Assume $n=4$. In this case, (\ref{ycw094})
can be equivalently written as
 \be\label{ycw095}
 (\cdots)(2+s)+12\alpha(b^2-4)^2\big[2(b^2-4)s^k_{\ 0|k}+64(4-b^2)\theta\alpha-2b^ks_{0|k}\big]=0,
 \ee
where the omitted term in the bracket is in a similar form as in
(\ref{ycw00}). By (\ref{ycw095}) we have
 \be\label{ycw096}
2(b^2-4)s^k_{\
0|k}+64(4-b^2)\theta\alpha-2b^ks_{0|k}=\xi\alpha(2+s),
 \ee
where $\xi=\xi(x)$ is a scalar function. Now it easily follows
from (\ref{ycw096}) that
 \be\label{ycw097}
s^k_{\ 0|k}=-16\theta\beta+\frac{b^ks_{0|k}}{b^2-4}.
 \ee
Using (\ref{ycw097}), (\ref{ycw095}) becomes
 \be\label{ycw098}
 Ric_{\alpha}=32\theta\alpha^2+\frac{4s_{0|0}}{b^2-4}-\frac{8(5b^2-8)s_0^2}{(b^2-4)^3}.
 \ee

Next we use Lemma \ref{lem0} to prove that $\beta$ is parallel
with respect to $\alpha$. Using (\ref{ycw084}), $t_{00}$ in
(\ref{ycw085}) and (\ref{ycw098}),  the fourth formula in Lemma
\ref{lem0} gives
 \be\label{ycw099}
 s^m_{|m}=\frac{16(2b^6+b^4-8b^2+32)\theta}{(b^2-4)(b^2+2)}.
 \ee
Using (\ref{ycw084}) and  $t_{00}$ in (\ref{ycw085}),  the third
formula in Lemma \ref{lem0} gives
 \be\label{ycw0100}
 b^ms_{0|m}=-8(b^2+2)\theta\beta.
 \ee
Using (\ref{ycw084}), $t_{00}$ in (\ref{ycw085}), (\ref{ycw098}),
(\ref{ycw099}) and  (\ref{ycw0100}), the second formula in Lemma
\ref{lem0} gives
 \be\label{ycw0101}
 s^m_{\ 0|m}=-\frac{24(b^4+4b^2-8)\theta}{(b^2-4)(b^2+2)}\beta.
 \ee
Plugging (\ref{ycw093}) and  (\ref{ycw0100}) into (\ref{ycw097})
yields
 \be\label{ycw0102}
 s^m_{\ 0|m}=-\frac{24(b^2-2)\theta}{b^2-4}\beta.
 \ee
Now (\ref{ycw0101}) and  (\ref{ycw0102}) imply $\theta=0$.
Therefore, by (\ref{ycw093}), we have $t^m_{\ m}=0$ and then
$\beta$ is closed. Further, by (\ref{ycw084}), $\beta$ is parallel
with respect to $\alpha$.

\

\noindent{\bf Case II:} Assume $n\ne 4$. In this case, by
(\ref{ycw094}) we have
 \be\label{ycw0103}
 \theta(b^2-4)^2\beta^2-s_0^2=c(4\alpha^2-\beta^2),
 \ee
 for some scalar function $c=c(x)$. By contractions on
 (\ref{ycw0103}), we have
  \be\label{ycw0104}
 c=-\theta b^2(b^2-4), \ \ \ (n-2)\theta b^2(b^2-4)=0.
  \ee

\noindent{\bf Case IIA:} Assume $n\ne2$. By (\ref{ycw0104}) we
have $\theta=0$. Then (\ref{ycw093}) shows $\beta$ is closed. So
in this case, $\beta$ is parallel with respect to $\alpha$ by
(\ref{ycw084}).

\

\noindent{\bf Case IIB:} Assume $n= 2$. Plugging (\ref{ycw0104})
into (\ref{ycw0103}) gives $s_0^2$ in (\ref{ycw085}). Since $n=2$,
we can always put
 \be\label{ycw0106}
 Ric_{\alpha}=\lambda\alpha^2,
 \ee
for some scalar function $\lambda=\lambda(x)$. Plug $s_0^2$ in
(\ref{ycw085}) and (\ref{ycw0106}) into (\ref{ycw094}) and then
(\ref{ycw094}) can be written as
 \be\label{ycw0107}
 (\cdots)\beta+8(b^2-4)\alpha^2\big[(b^2-4)s^k_{\
 0|k}-b^ks_{0|k}\big].
 \ee
So by (\ref{ycw0107}) we have
 \be\label{ycw108}
 s^k_{\ 0|k}=\tau\beta+\frac{b^ks_{0|k}}{b^2-4},
 \ee
where $\tau=\tau(x)$ is a scalar function. Plug (\ref{ycw108})
into (\ref{ycw0107}) and then (\ref{ycw0107}) becomes
 \be\label{ycw109}
 s_{0|0}=\Big\{\frac{b^2-4}{3}(\lambda-\tau)+\frac{8(13b^4+4b^2-80)\theta}{3(b^2-4)}\Big\}\alpha^2-
 \Big\{\frac{b^2-4}{12}(\lambda+2\tau)+\frac{8(b^4+19b^2-56)\theta}{3(b^2-4)}\Big\}\beta^2.
 \ee
Using the contraction $b^ib^js_{i|j}$, it follows from
(\ref{ycw109}) that
 $$
 \tau=-\frac{16(b^4+7b^2-26)\theta}{(2+b^2)(b^2-4)}-\frac{(b^2-4)\lambda}{2(2+b^2)}.
 $$
 Plugging the above into (\ref{ycw108}) gives (\ref{ycw086}) and into (\ref{ycw109}) gives (\ref{ycw087}).
 Finally, using (\ref{ycw084})-(\ref{ycw087}), the flag curvature
 ${\bf K}$ of $F$ is isotropic and it is given by (\ref{ycw088}).
                                                    \qed

\section{General case}

We first give the condition in the following proposition for a
$p$-power metric $F=\alpha(1+\beta/\alpha)^p$ to be positively
definite.

\begin{prop}\label{prop71}
A $p$-power metric $F=\alpha(1+\beta/\alpha)^p$ ($p\ne 0$) is
positively definite if and only if one of the following cases
holds:
 \ben
  \item[{\rm (i)}] ($p> 2$, or $p<0$) $b^2$ satisfies
      $$b^2<\frac{1}{(p-1)^2}.$$
  \item[{\rm (ii)}] ($\frac{1}{2} \le p\le 2$) $b^2$ satisfies
         $$b^2<1.$$

   \item[{\rm (iii)}] ($0 < p< \frac{1}{2}$) $b^2$ satisfies
    $$
        b^2< \frac{(2-p)^2}{4(1-p^2)^2}.
    $$
 \een
\end{prop}

{\it Proof :}  It's an elementary discussion. We only need to
start from the following conditions
 $$
 \phi(s)>0,\ \ \ \ \   \phi(s)-s\phi'(s)>0, \ \ \ \ \
 \phi(s)-s\phi'(s)+(b^2-s^2)\phi''(s)>0,
 $$
where $\phi(s):=(1+s)^p$, $s:=\beta/\alpha$ and
$b:=||\beta||_{\alpha}$.   Details are omitted.          \qed

\

In the following, we will prove Theorem \ref{th1} when $F$ is a
$p$-power metric $F=\alpha(1+\beta/\alpha)^p$ with
 \be\label{ycw00109}
 p\ne 1,2,-1,\frac{1}{2}.
 \ee

 \begin{prop}
An $n$-dimensional $p$-power metric $F=\alpha(1+\beta/\alpha)^p$
satisfying (\ref{ycw00109}) is Einstein-reversible if and only if
$F$ is Ricci-flat-parallel.
\end{prop}

{\it Proof :} For the $p$-power metric
$F=\alpha(1+\beta/\alpha)^p$, let
 $f(x,y)$ defined by (\ref{cw6}).
Then by the definition of Einstein-reversibility, we have
$f(x,y)=f(x,-y)$, which multiplied by
 $$\alpha^4[1-(p-1)^2s^2]^3\big\{[1+p(p-1)b^2-(p^2-1)s^2]^2-(p-2)^2s^2]^4\big\}$$
can be firstly written in the following form
 \beq\label{ycw0108}
 0&=&\Big[(\cdots)+(\cdots)\big(\frac{1+s}{1-s}\big)^{2p}\Big]\big[1+p(p-1)b^2-(p-2)s-(p^2-1)s^2\big]
 +\nonumber\\
 &&(p-2)^4s^4(1+s)^2[1-(p-1)s]^3[1+(p-1)s]^3[2(p^2-1)s+p-2]^2\times\nonumber\\
 &&\big[(1-ps+s)r_{00}-2p\alpha s_0\big]^2,
 \eeq
where  the omitted terms in the brackets are in a similar form as
in (\ref{ycw00}). We discuss (\ref{ycw0108}) in the following two
cases, and we only need to assume $2p$ is an integer.

\

\noindent{\bf Case I:} Assume the following holds
 \be\label{ycw0109}
 \alpha\big[(1-ps+s)r_{00}-2p\alpha
 s_0\big]=(g_0+g_1\alpha)\alpha^2\big[1+p(p-1)b^2-(p-2)s-(p^2-1)s^2\big],
 \ee
where $g_0$ is a 1-form and $g_1=g_1(x)$ is a scalar function. By
(\ref{ycw0109}), we first get
 \be\label{ycw0110}
 r_{00}=g_1\big[(1+p^2b^2-pb^2)\alpha^2-(p^2-1)\beta^2\big]-(p-2)g_0\beta,\
 \ \ g_0=-\frac{(p^2-1)g_1}{2p-1}\beta.
 \ee
Plugging (\ref{ycw0110}) into (\ref{ycw0109}) gives
 \be\label{ycw0111}
 (p-2)g_1\big[(p-1)^2b^2-1\big]\beta+2(2p-1)s_0=0.
 \ee
Using $b^ms_m=0$, it follows from (\ref{ycw0111}) that $s_0=0$ and
$g_1=0$. So by (\ref{ycw0110}) we have
 \be\label{ycw0112}
 r_{00}=0, \ \ \ \  s_0=0.
 \ee
Now using (\ref{ycw0112}), (\ref{ycw0108}) can be equivalently
written as
 \be\label{ycw0113}
 \Big[(\cdots)+(\cdots)\big(\frac{1+s}{1-s}\big)^{2p}\Big]\big[1-(p-1)s\big]+p^2t_{00}=0,
 \ee
where  the omitted terms in the brackets are in a similar form as
in (\ref{ycw00}). By (\ref{ycw0113}) we have
 \be\label{ycw0114}
 t_{00}=(g_0+\theta\alpha)\alpha\big[1-(p-1)s\big],
 \ee
where $g_0$ is a 1-form and $\theta=\theta(x)$ is a scalar
function. It follows from (\ref{ycw0114}) that
 \be\label{ycw0115}
 t_{00}=\theta\big[\alpha^2-(p-1)^2\beta^2\big].
 \ee
Since $s_0=0$ by (\ref{ycw0112}), we have $t_0=0$, and then a
contraction on (\ref{ycw0115}) shows $\theta=0$. Thus $\beta$ is
closed. Plus (\ref{ycw0112}), it shows that $\beta$ is parallel
with respect to $\alpha$.

\

\noindent{\bf Case II:} Assume $1+p(p-1)b^2-(p-2)s-(p^2-1)s^2$ has
one of the following factors
  $$
 s,\ \  1\pm s,\ \  1\pm (p-1)s, \ \ 2(p^2-1)s+p-2.
  $$
In this case, $b^2$ must be a constant. By the positively definite
condition shown in Proposition \ref{prop71}, we only need to
consider three cases:
 $$
 b^2=\frac{4-p}{p(p-1)^2},\ \ \ b^2=\frac{4-5p}{4(p-1)^2(p+1)},\
 \ \ b^2=\frac{p^2+p-4}{p(p-1)}.
 $$
We only investigate one case, and the left is similar. Suppose
\be\label{ycw0116}
 b^2=\frac{4-p}{p(p-1)^2}.
 \ee
  Then
 $$
 1+p(p-1)b^2-(p-2)s-(p^2-1)s^2=-\frac{\big[1+(p-1)s\big](p^2s-s-3)}{p-1}.
 $$
So in (\ref{ycw0108}), $p^2s-s-3$ must be a factor of
$(1-ps+s)r_{00}-2p\alpha s_0$, which implies
 \be\label{ycw0117}
 (1-ps+s)r_{00}-2p\alpha
 s_0=(g_0+g_1\alpha+g_2\alpha^2)\alpha(p^2s-s-3),
 \ee
where $g_0$ is a polynomial of degree two, $g_1$ is a 1-form and
$g_2=g_2(x)$ is a scalar function. By (\ref{ycw0117}) we can
easily get
 \be\label{ycw0118}
 r_{00}=\frac{2p(p-1)(p+1)^2}{3(p-2)}\beta
 s_0+\frac{9\alpha^2-(p^2-1)^2\beta^2}{3(p-1)(p-2)}\ c,
 \ee
where $c=c(x)$ is a scalar function. Since $b^2=constant$, we have
$r_0+s_0=0$. By a contraction on (\ref{ycw0118}) we can get $r_0$,
and then using (\ref{ycw0116}) and $r_0+s_0=0$ we obtain
 $$c(p-2)\beta-p(p-5)s_0=0.$$
By (\ref{ycw0116}), $p\ne 5$. Using $b^ms_m=0$, the above shows
$c=0,\ s_0=0$. Then by (\ref{ycw0118}), we again get
(\ref{ycw0112}). Thus by  a similar proof as that in Case I, we
have $t_{00}=0$. So $\beta$ is parallel with respect to $\alpha$.
\qed

\section{Proof of Theorem \ref{th2}}

\begin{thm}\label{th3}
 Let $F=\sqrt{\alpha(\alpha+\beta)}$ be a two-dimensional
 square-root metric. Then $F$ is Einsteinian (equivalently, of
 isotropic flag curvature) if and only if (\ref{ycw084}) holds. In
 this case, the flag curvature of $F$ is given by
  \be\label{ycw0123}
 {\bf
 K}=\frac{2}{2+b^2}\Big\{\lambda-\frac{8s_ms^m}{b^2(b^2-4)}\Big\},
  \ee
  where $\lambda$ is the sectional curvature of $\alpha$.
\end{thm}

{\it Proof :} Assume (\ref{ycw084}) holds. By Proposition
\ref{prop61}, we only need to show that
(\ref{ycw085})-(\ref{ycw087}) automatically hold.

Since $n=2$, we can always put (\ref{ycw0106}).

Fix a point $x\in M$ and take  an orthonormal basis
  $\{e_i\}$ at $x$ such that
   $$\alpha=\sqrt{(y^1)^2+(y^2)^2},\ \ \beta=by^1.$$
Then
$$s_0^2=s_2y^2,\ \ \ b^2\alpha^2-\beta^2=b^2y^2.$$
So for some $\theta=\theta(x)$, the second formula in
(\ref{ycw085}) holds. Then we immediately get
 \be\label{ycw0124}
 s_ms^m=4b^2(b^2-4)\theta.
 \ee

Since $n=2$, we always have
 \be\label{ycw0125}
 s_{ij}=\frac{b_is_j-b_js_i}{b^2}.
 \ee
Therefore, by (\ref{ycw0125}), (\ref{ycw0124}) and the second
formula in (\ref{ycw085}), we have
 $$t_{00}=-\frac{s_ms^m\beta^2+b^2s_0^2}{2}=-4(b^2-4)\theta\alpha^2,$$
which is just the first formula in (\ref{ycw085}).

Using (\ref{ycw084}), (\ref{ycw093}) and (\ref{ycw0106}), it
follows from the fourth formula in Lemma \ref{lem0} that
 \be\label{ycw0126}
 s^m_{|m}=\frac{32(b^2-1)(b^4+4b^2-8)\theta}{(b^2-4)(b^2+2)}+\frac{b^2(b^2-4)\lambda}{2(b^2+2)}.
 \ee
Similarly, using (\ref{ycw084}) and  (\ref{ycw093}), it follows
from the third formula in Lemma \ref{lem0} that
 \be\label{ycw0127}
 b^ms_{0|m}=-8(b^2+2)\theta\beta.
 \ee
Now using (\ref{ycw0126}) and (\ref{ycw0127}), plugging
(\ref{ycw084}), (\ref{ycw093}) and (\ref{ycw0106}) into the second
formula in Lemma \ref{lem0} gives
 \be\label{ycw0128}
 s^m_{\
 0|m}=-\Big\{\frac{24(b^2+8)(b^2-2)\theta}{(b^2-4)(b^2+2)}+\frac{(b^2-4)\lambda}{2(b^2+2)}\Big\}\beta.
 \ee
Plug (\ref{ycw0127}) into (\ref{ycw086}) and then (\ref{ycw086})
is just (\ref{ycw0128}).

Finally, we prove (\ref{ycw087}). Since $\alpha$ is locally
conformally flat, $\alpha$ and $\beta$ can be locally expressed as
 \be\label{ycw0129}
 \alpha=e^{\sigma}\sqrt{(y^1)^2+(y^2)^2},\  \ \
 \beta=e^{\sigma}(\xi y^1+\eta y^2),
 \ee
 where $\sigma=\sigma(x),\ \xi=\xi(x),\ \eta=\eta(x)$ are scalar
 functions. In the following, we define
 $$
\sigma_{i}:=\sigma_{x^i},\ \ \ \sigma_{ij}:=\sigma_{x^ix^j}, \ \
etc.
 $$
The sectional curvature $\lambda$ of $\alpha$ and the
 norm $b=||\beta||_{\alpha}$ are given by
  \be\label{ycw0130}
 \lambda=-e^{-2\sigma}(\sigma_{11}+\sigma_{22}),\ \ \
 b^2=\xi^2+\eta^2.
  \ee
   Further, we have
   \be\label{ycw0131}
  s_ms^m=\frac{(\xi^2+\eta^2)(\xi_2+\xi\sigma_2-\eta\sigma_1-\eta_1)^2}{4e^{2\sigma}}.
   \ee
 Without loss of generality we may assume $\xi\ne 0$. A
 direct computation shows (\ref{ycw084}) is equivalent to
 \beq
&&\xi_1=-\frac{\eta(\xi\xi_2+\eta\eta_2+\xi\eta_1)}{\xi^2},\label{ycw0132}\\
 &&\sigma_1=\frac{3\xi\eta\xi_2-\eta_2(2\xi^2-\eta^2-2)}{2\xi(\xi^2+\eta^2-1)},\
 \ \ \
 \sigma_2=\frac{(2\eta^2-\xi^2-2)(\xi\xi_2+\eta\eta_2)}
 {2\xi^2(\xi^2+\eta^2-1)}+\frac{\eta_1}{\xi}.\label{ycw0133}
 \eeq
 For the system (\ref{ycw0132}) and (\ref{ycw0133}), its
 integrable condition $\sigma_{12}=\sigma_{21}$ is given by
  \beq\label{ycw0134}
 0&=&2\eta\xi^2\xi_2^2+\big[3\eta\xi^2\eta_1+\xi(5\eta^2-2\xi^2)\eta_2\big]\xi_2
 +\eta\xi^2\eta_1^2+3\xi\eta^2\eta_1\eta_2+\nonumber\\
 &&2\eta(\eta^2-\xi^2)\eta_2^2-\xi^2(\xi\eta\xi_{22}+\eta^2\eta_{22}-\xi^2\eta_{22}-\xi^2\eta_{11}).
  \eeq
Now plug (\ref{ycw0129})-(\ref{ycw0131}) into (\ref{ycw087}) and
then we write (\ref{ycw087}) in the following form
 $$\widehat{A}(x)(y^1)^2+\widehat{B}(x)y^1y^2+\widehat{C}(x)(y^2)^2=0.$$
 So (\ref{ycw087}) holds if and only if
  $$\widehat{A}(x)=0,\ \ \ \widehat{B}(x)=0,\ \ \ \widehat{C}(x)=0.$$
  Plugging (\ref{ycw0132}) and (\ref{ycw0133}) into $\widehat{A}(x),\  \widehat{B}(x)$
  and
  $\widehat{C}(x)$, a direct computation shows that (i) $\widehat{A}(x)=0\Leftrightarrow
  \widehat{B}(x)=0$; (ii) $\widehat{C}(x)=0$; and further (iii) $\widehat{A}(x)=0\Leftrightarrow
  (\ref{ycw0134})$. Therefore, (\ref{ycw084}) implies
  (\ref{ycw087}).

  Eq. (\ref{ycw0123}) directly follows from (\ref{ycw088}) and
  (\ref{ycw0124}).      \qed

  \

  \noindent {\it Proof of Theorem \ref{th2} :}

  \

 We first take a deformation on $\beta$ but keep $\alpha$
unchanged, and then we get a Killing vector field. Define a
Riemannian metric $\widetilde{\alpha}$ and 1-form
$\widetilde{\beta}$ by
 \be\label{ycw0135}
\widetilde{\alpha}:=\alpha, \ \
\widetilde{\beta}:=(1-b^2)^{-\frac{3}{4}}\beta.
 \ee
By (\ref{ycw0135}), we show that (\ref{ycw084}) is reduced to
 \be\label{ycw0136}
 \widetilde{r}_{ij}=0.
 \ee
So $\widetilde{\beta}$ is a Killing form.

Now we  express $\alpha$ locally as
 \be\label{ycw0137}
\alpha:=e^{\sigma}\sqrt{(y^1)^2+(y^2)^2},
 \ee
where $\sigma=\sigma(x)$ is a scalar function.  Then by the result
in \cite{Y}, we have
 \be\label{ycw0138}
 \widetilde{\beta}=\widetilde{b}_1y^1+\widetilde{b}_2y^2=e^{2\sigma}(uy^1+vy^2),
 \ee
where $u=u(x),v=v(x)$ are a pair of scalar functions such that
 $$f(z)=u+iv, \ \ z=x^1+ix^2$$
 is a complex analytic function, and further by (\ref{ycw0136}),
 $u,\ v$ and $\sigma$ satisfy the following PDEs:
  \be\label{ycw0139}
 u_1=v_2,\ \ \ \ u_2=-v_1,\ \ \ \ u_1+u\sigma_1+v\sigma_2=0.
  \ee

 Next we determine $\sigma$ in terms of the triple $(B,u,v)$, where $B:=b^2$.
 Firstly by (\ref{ycw0135}) and then by (\ref{ycw0137}) and
(\ref{ycw0138}) we get
 \be\label{ycw0140}
||\widetilde{\beta}||^2_{\alpha}=\frac{B}{(1-B)^{\frac{3}{2}}},\ \
||\widetilde{\beta}||^2_{\alpha}=e^{2\sigma}(u^2+v^2).
 \ee
Therefore, by (\ref{ycw0140}) we get
 \be\label{ycw0141}
 e^{2\sigma}=\frac{B}{(u^2+v^2)(1-B)^{\frac{3}{2}}}.
 \ee
Now it follows from (\ref{ycw0135}), (\ref{ycw0137}),
(\ref{ycw0138}) and (\ref{ycw0141}) that (\ref{ygj001}) holds.

Comparing (\ref{ycw0129}) with (\ref{ygj001}), we have
 \be\label{ycw0142}
 \xi=\frac{u\sqrt{B}}{\sqrt{u^2+v^2}},\ \ \ \eta=\frac{v\sqrt{B}}{\sqrt{u^2+v^2}},
 \ \ \
 \sigma=ln\frac{\sqrt{B}}{\sqrt{u^2+v^2}(1-B)^{\frac{3}{4}}}.
 \ee
Using $\sigma$ in (\ref{ycw0142}) and $u_1=v_2,\ u_2=-v_1$ in
(\ref{ycw0139}), $u_1+u\sigma_1+v\sigma_2=0$ is reduced to
 $$uB_1+vB_2=0.$$
So by (\ref{ycw0139}) we obtain (\ref{ygj002}).

 Finally, we prove
that the flag curvature ${\bf K}$ is given by (\ref{ygj003}).
Plugging  (\ref{ycw0142}) into (\ref{ycw0131}) and using
$u_1=v_2,\ u_2=-v_1$, we have
 \be\label{ycw0143}
 s_ms^m=\frac{(B-4)^2(uB_2-vB_1)^2}{64B\sqrt{1-B}}.
 \ee
By (\ref{ygj002}) we have
 \be\label{ycw0144}
 u_1=v_2,\ \ \ u_2=-v_1,\ \ \ u_{11}+u_{22}=0,\ \ \ v_{11}+v_{22}=0,\ \ \ B_2=-\frac{u}{v}B_1.
 \ee
Plugging  (\ref{ycw0141}) into (\ref{ycw0130}) and using
(\ref{ycw0144}), we get
 \be\label{ycw0145}
\lambda=-\frac{u^2+v^2}{4B^2}\Big\{(B+2)\sqrt{1-B}(B_{11}+B_{22})+\frac{(u^2+v^2)(B^2+4B-2)B_1^2}{Bv^2\sqrt{1-B}}\Big\}.
 \ee
Now plugging (\ref{ycw0143}) and (\ref{ycw0145}) into
(\ref{ycw0123}) and using (\ref{ycw0144}), we obtain
(\ref{ygj003}).      \qed

\vspace{0.6cm}

\noindent Guojun Yang \\
Department of Mathematics \\
Sichuan University \\
Chengdu 610064, P. R. China \\
yangguojun@scu.edu.cn

\end{document}